\documentclass{amsart}

\usepackage{amsfonts}

\usepackage{graphics}

\font\sr = cmr8
\font\br = cmr12

\newtheorem{theorem}{Theorem}

\theoremstyle{plain}

\newtheorem{claim}{Claim}

\newtheorem{corollary}{Corollary}

\begin{document}

\title{Incompressible surfaces and $(1,1)$-knots}

\author {Mario Eudave-Mu\~noz}

\address {Instituto de Matem\'aticas, UNAM,  Circuito Exterior,
Ciudad Universitaria, \break  04510 M\'exico D.F., MEXICO.
E-mail: {\rm mario@matem.unam.mx}}


\begin{abstract}
Let $M$ be $S^3$, $S^1\times S^2$, or a lens space $L(p,q)$, and let $k$ be a
$(1,1)$-knot in $M$. We show that if there is a closed meridionally
incompressible surface in the complement of  $k$, then the surface and the knot can
be put in a special position, namely, the surface is the boundary of a regular
neighborhood of a toroidal graph, and the knot is level with respect to that graph.
As an application we show that for any such $M$ there exist tunnel number one knots
which are not $(1,1)$-knots.
\end{abstract}

\maketitle

\section{Introduction}

An important problem in knot theory is that of determining all 
incompressible surfaces in a given knot complement.  
The purpose of this paper is to give a characterization of all $(1,1)$-knots
which admit a closed meridionally incompressible surface in their complement.
The family of $(1,1)$-knots is an interesting class of knots
which has received considerable attention recently 
(see for example \cite{CM}, \cite{CK}, \cite{GMM}). Recall that $k$ is a 
$(1,1)$-knot if there
is a standard torus $T$ in $S^3$ such that $k$ is of 1-bridge with respect to $T$, 
that is, $k$ intersects $T$ in two points, which divide $k$ into two arcs which are
parallel to arcs lying on $T$. The family of $(1,1)$-knots contains the torus knots,
2-bridge knots, and is contained in the
class of tunnel number one knots. Recall that a  knot $k$ has
tunnel number one if there exists an arc $\tau$ embedded in
$S^3$, with $k\cap \tau=\partial \tau$, such that $M-int\, N(k\cup \tau)$
is a genus 2-handlebody.

Some results are known about incompressible surfaces in complements of
tunnel number one knots, and in particular about $(1,1)$-knots.
It is known that the complement of a tunnel number one knot
does not contain any essential meridional planar surface, i.e
a planar surface which is incompressible, $\partial$-incompressible
and whose boundary consists of meridians of the knot \cite{GR}.
There are tunnel number one knots whose complement contains an
essential torus; such knots
have been classified by K. Morimoto and M. Sakuma \cite{MS}
(see also \cite{E1}), and it turns out that
all of them are $(1,1)$-knots. The author constructed for any $g \geq 2$, infinitely
many tunnel number one knots whose complement contains a closed incompressible
and meridionally incompressible surface of genus $g$ \cite{E2}. Many of these
 knots are
$(1,1)$-knots.

The construction given in \cite{E2} is roughly as follows (see Section 2 for more
details):

Take $n$ parallel tori $T_1,\dots,T_n$, standardly embedded
in $S^3$, and take an essential curve $\gamma_i$ in each one of
the tori. Take $n-1$ straight arcs $\alpha_i$ joining consecutive curves.
The union of the curves and the arcs form a graph $\Gamma$, which
we call a toroidal graph. Note that $N(\Gamma)$ is a handlebody of
genus $n$. If the curves $\gamma_i$ satisfy certain natural conditions 
then $S=\partial N(\Gamma)$ is 
incompressible in $S^3-int\, N(\Gamma)$.
Let $\eta(\gamma_i)$ be a neighborhood of $\gamma_i$ in the torus
$T_i$, this is an annulus. Now let $k$ be a knot in level position with
respect to $\Gamma$; by this we mean that
$k\subset (\cup \eta(\gamma_i))\cup  (\cup N(\alpha_i))$,
such that $k\cap \eta(\gamma_i)$ consists of two arcs for $2\leq i\leq n-1$,
$k\cap \eta(\gamma_i)$ is one arc for $i=1,n$, and $k \cap N(\alpha_i)$
consists of two straight arcs. It follows that $k$ is a $(1,1)$-knot,
and if $k$ is well wrapped in $N(\Gamma)$, then $S$ is
incompressible in $N(\Gamma)-k$.

In Section 3 we show that if $k$ is a $(1,1)$-knot and $S$ is a closed
meridionally incompressible surface in its complement, then
$k$ and $S$ come from the above construction. To do that, consider first the
following equivalent definition of $(1,1)$-knots. Let $T\times I\subset S^3$, 
where $I=[0,1]$, and $T$ is an standard torus. A knot $k$ is a $(1,1)$-knot
if and only if $k$ lies in $T\times I$, so that $k\cap (T\times \{ 0\}) = k_0$ 
is an arc,
$k\cap (T\times \{ 1\}) = k_1$ is an arc, and $k\cap (T\times (0,1))$ consists
of two straight arcs, i.e., arcs which cross transversely every torus
$T\times \{x\}$. Let $k$ be a knot in such a position and let $S$ be a
closed meridionally incompressible surface in the complement of
$k$. $S$ can be isotoped, keeping $k$ fixed, so that most of $S$ lie in the product
$T\times I$. Now put $S$ in Morse position in the product $T\times I$. 
By analyzing the singular points of the surface we show that it can be positioned
so that it looks like $\partial N(\Gamma)$ for a certain toroidal graph $\Gamma$,
and that $k$ is in level position with respect to $\Gamma$. 

The class of $(1,1)$-knots can be defined for any manifold $M$ having a Heegaard
decomposition of genus 1, that is, it can be defined for 
$M=S^3, S^1\times S^2$ or a
lens space $L(p,q)$. The definition is the same, where $T$ is supposed to be a
Heegaard torus. Our results are also valid for $(1,1)$-knots in these manifolds.

Our results imply that if in the complement of a
$(1,1)$-knot $k$ in $M$ there is a closed meridionally
incompressible surface, then the surface bounds a handlebody in $M$ in which $k$
lies. In general the knots constructed in Section 6 of \cite{E2} do not have this
property. So we can construct for any such $M$, many tunnel number one knots which
are  not $(1,1)$-knots. See Corollary 1 in Section 2.

\section{Statement of results}

Let $M$ be a compact orientable 3-manifold, and let $S \not= S^2$ be a 
compact orientable surface in $M$, either properly embedded or contained in 
$\partial M$.
Let $k$ be a knot in $M$ disjoint from $S$. We say that $S$ is meridionally 
compressible  
if there is an embedded disk $D$ in $M$, with $S\cap D=\partial D$, 
which is a nontrivial curve in $S$, and so that $k$ intersects $D$ 
at most in one point. Otherwise $S$ is called meridionally incompressible.
In particular, if $S$ is meridionally incompressible then it is
incompressible in $M-k$.
If the surface $S$ is a 2-sphere, then $S$ is usually defined to be incompressible
if it is not the boundary of a 3-ball; similarly we define it to be meridionally
incompressible if it is incompressible.

In what follows $M$ is $S^3, S^1\times S^2$, or a lens space $L(p,q)$.
If $k$ is a $(1,1)$-knot in $M$ (or in general, a tunnel number one knot), 
and $S$ is an incompressible sphere in $M-k$, then it
can be proved that $k$ is a trivial knot, i.e., $k$ is the boundary of a disk. So
there are two cases of $(1,1)$-knots which contain an incompressible sphere: 1)
$M=S^1\times S^2$, or $L(p,q)$, $k$ a trivial knot, and $S$ a sphere in $M$ bounding
a ball in which $k$ lies. 2) $M=S^1\times S^2$, $k$ a trivial knot, and $S$ a
non-separating sphere in $M-k$. From now on we will assume that we are not in this
situation, that is, we assume that all our surfaces are of genus $\geq 1$.
If $S$ is a closed non-separating surface in the complement of a tunnel 
number one 
knot $k$ in $M$, then by doing an argument similar to that in \cite{J}, 
\cite{S}, or \cite{E1}, and  applied to our situation, it follows that 
either the knot is trivial or there is a closed non-separating surface 
in the complement of $k$ disjoint from an unknotting tunnel, which is 
clearly impossible. So assume in what follows 
that all surfaces are separating.

In \cite{E2} some constructions are given of tunnel number one knots whose 
complements
contain a closed meridionally incompressible surface. The
simplest of those knots are also $(1,1)$-knots. We describe this
construction. 

Let $T$ be a Heegaard torus in $M$, and let $I=[0,1]$. Consider 
$T\times I\subset M$.
$T\times \{ 0\}=T_0$ bounds a solid torus $R_0$, and $T_1=T\times \{1\}$ bounds a
solid torus $R_1$, such that $M= R_0\cup (T\times I)\cup R_1$. Choose $n$ distinct
points on $I$, $e_1=0,\ e_2,\dots, e_n=1$, so that $e_i < e_{i+1}$, for all 
$i\leq n-1$. Consider the tori
$T\times \{e_i\}$ (so $T\times \{e_1\}= T_0$, $T\times \{e_n\}= T_1$). There is a
natural projection of each of these tori into $T\times \{e_1\}$. 
A slope in $T$ is a class of isotopy of essential simple closed curves. Slopes in
$T\times \{x\}$ are in one to one correspondence with slopes in $T\times \{y\}$, for
any $x,y\in I$, so comparing curves on different tori is accomplished via this 
isotopy.
If $\beta$ and $\gamma$ are two slopes on $T$, $\Delta(\beta,\gamma)$ denotes, as
usual, their minimal geometric intersection number. 
Denote by $\mu_M$ the essential simple closed curve on 
$T\times \{e_n\}$ which bounds a meridional disk in $R_1$, and denote by $\lambda_M$
the essential simple closed curve on $T\times \{e_1\}$ which bounds a disk
in $R_0$. So $\Delta(\mu_M,\lambda_M)=1$ if
$M=S^3$, $\Delta(\mu_M,\lambda_M)=0$ if $M=S^1\times S^2$, and
$\Delta(\mu_M,\lambda_M)=p$ if $M=L(p,q)$. By a straight arc in a product 
$T\times [a,b]$ we mean a properly embedded arc which is transverse to every torus 
$T\times \{x\}$  in the product.  

Let $\gamma_i$ be a simple closed essential curve embedded in the
torus $T\times \{e_i\}$.
Let $\alpha_i$, for $1\leq i \leq n-1$, be a 
straight arc in $T\times [e_i,e_{i+1}]$, joining $\gamma_i$
and $\gamma_{i+1}$, and assume that $\alpha_i \cap \alpha_j=\emptyset$
if $i\not= j$.
Let $\Gamma$ be the 1-complex consisting of the union of all the
curves $\gamma_i$ and the arcs $\alpha_j$. So $\Gamma$ is a trivalent
graph embedded in $M$. We say that $\Gamma$ is a toroidal graph of type
$n$, $n\geq 2$,  if it satisfies the following:

\begin{enumerate}
\item $\Delta(\gamma_1,\lambda_M) \geq 2$, that is, $\gamma_1$ is not
homotopic to the core of the solid torus $R_0$.

\item $\Delta(\gamma_i,\gamma_{i+1}) \geq 2$, for all $i\leq n-1$.

\item $\Delta(\gamma_n,\mu_M) \geq 2$, that is, $\gamma_n$ is not
homotopic to the core of the solid torus $R_1$.

\end{enumerate}

We define $\Gamma$ to be a toroidal graph of type 1, if it consists of a 
single curve
$\gamma_1$ in $T\times \{e_1\}$, such that $\Delta(\gamma_1,\lambda_M) \geq 2$
and $\Delta(\gamma_1,\mu_M) \geq 2$. See
Figure 1 for an example of a toroidal graph of type 2 in $S^3$.

\bigskip

\centerline{\includegraphics{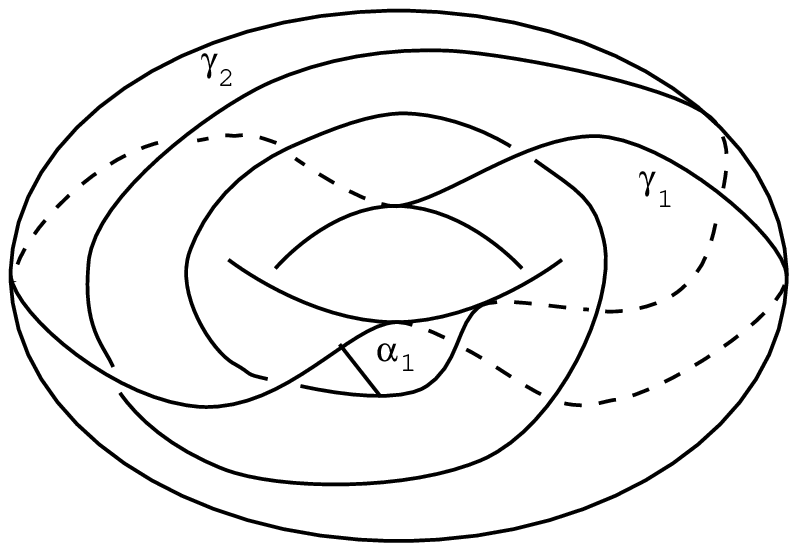}}

\centerline{\sr Figure 1}

\bigskip

Let $N(\Gamma)$ be a regular neighborhood of $\Gamma$. This is
a genus $n$ handlebody. We can assume that $N(\Gamma)$ is made
of the union of $n$ solid tori $N(\gamma_i)$, which are joined
by $(n-1)$ 1-handles $N(\alpha_j)$. Note that there are disks
$D_1,\dots,D_{n-1}$, properly embedded in $N(\Gamma)$, with $D_j$
contained in $N(\alpha_j)$ and intersecting $\alpha_j$ once. These disks
cut $N(\Gamma)$ into $n$ solid tori.

\begin{theorem} Let $\Gamma$ be a toroidal graph in $M$. Then
$S=\partial N(\Gamma)$ is incompressible in $M-int\,N(\Gamma)$.
\end{theorem}

This is proved in Theorem 4.1 of \cite{E2} for the case $M=S^3$; the proof in the 
general case is
essentially the same. Note that if any of the conditions for a toroidal graph fails,
then $S$ will be compressible.

Let $\Gamma= (\cup \gamma_i)\cup (\cup \alpha_i)$ be a toroidal graph
of type $n$, with $n\geq 2$.
Denote by $\eta(\gamma_i)$ a regular neighborhood
of $\gamma_i$ in the torus $T\times \{e_i\}$, this is an annulus, we can assume that
$N(\gamma_i)\cap (T\times \{e_i\})=\eta(\gamma_i)$.

A knot $k$ contained in $T\times I$ is level with respect to $\Gamma$
if it satisfies the following:

\begin{enumerate}
\item $k\cap (T\times \{e_1\})$ consists of one arc contained in $int\,\eta(\gamma_1)$.

\item $k\cap (T\times \{e_i\})$ consists of two arcs contained in $int\,\eta(\gamma_i)$,
for $2\leq i\leq n-1$.

\item $k\cap (T\times \{e_n\})$ consists of one arc contained in $int\,\eta(\gamma_n)$.

\item $k\cap (T\times (e_i,e_{i+1}))$ consists of two straight arcs
contained in $int\,N(\alpha_i)$, for $1\leq i\leq n-1$.

\end{enumerate}

In the special case of a toroidal graph of type 1, we say that $k$ is level
with respect to $\Gamma$ if: (1) $k\cap T_0$ is a single arc contained in
$\eta (\gamma_1)$,
and (2) $k\cap (T\times (0,1])$ is a single arc which can be projected to
an arc on $\eta (\gamma_1)$. So if a knot $k$ is level with respect to a
toroidal graph of
type 1, it is either a torus knot, or a satellite of a torus knot.

It is not difficult to see that if $k$ is level with respect to a toroidal
graph $\Gamma$ then $k$ is a $(1,1)$-knot; to see this, just isotope $k$ so that
$k\cap (T\times (0,1))$ consists of two straight arcs.
By construction, if $k$ is level with respect to a toroidal graph $\Gamma$,
then $k$ lies in $N(\Gamma)$, and $k$ intersects each disk $D_i$
twice. By cutting $k$ along the disks $D_i$, and then joining its endpoints
with an arc on $D_i$, we get $n$ knots $k_1,\dots, k_n$ embedded in
$N(\alpha_1),\dots,N(\alpha_n)$ respectively. We say that $k$ is well wrapped in
$N(\Gamma)$ if the wrapping number of each $k_i$ in $N(\alpha_i)$ is $\geq 2$.

See Figure 2 for an example of a knot which is level with respect to the
toroidal graph in Figure 1.

\medskip

\centerline{\includegraphics{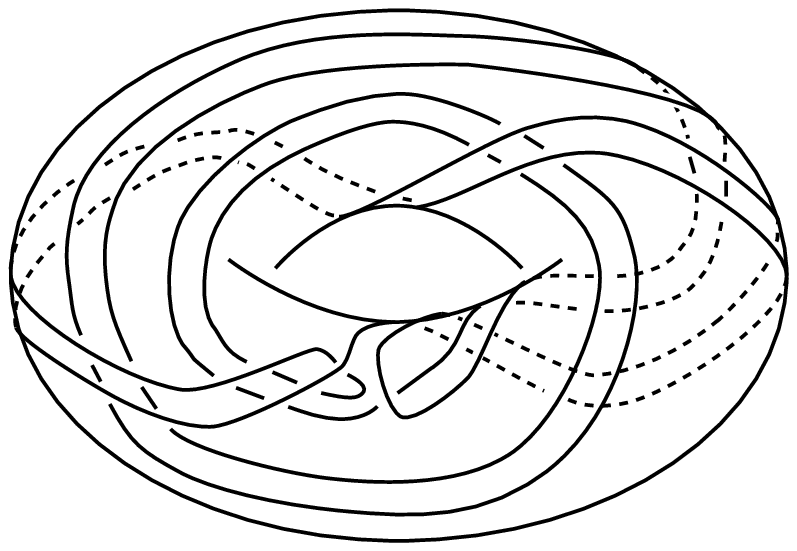}}

\centerline{\sr Figure 2}

\bigskip

\begin{theorem} Let $k$ be a knot which is level with respect to a toroidal
graph $\Gamma$ of type $n$. Assume also that $k$ is well wrapped in $N(\Gamma)$.
Then the closed surface $S=\partial N(\Gamma)$ is 
meridionally incompressible in $M-k$.
\end{theorem}

This is proved in Theorem 5.1 of \cite{E2}. It is not difficult to see that if $k$ 
is level with respect to $\Gamma$, but it is not well wrapped in $N(\Gamma)$ then 
$S$ is compressible in $N(\Gamma)-k$.

In this paper we prove the following,

\begin{theorem} Let $k$ be a $(1,1)$-knot in $M$. Let $S$ be a closed
meridionally incompressible surface in $M-k$. Then there is a toroidal 
graph $\Gamma$,
so that $S=\partial N(\Gamma)$, and $k$ is level with respect to $\Gamma$.
\end{theorem}

In the special case that the surface $S$ is a torus, this result has been
previously proved in \cite{H}, where all the $(1,1)$-knots which contain
an incompressible torus are determined. In the case that $M=S^3$, this result
follows from \cite{MS} and \cite{E1}, for in those papers all tunnel number 
one knots in
$S^3$ which contain an essential torus are determined, and all of them turn out to
be $(1,1)$-knots (see 4.8 in \cite{E1}).

An important problem about $(1,1)$-knots and meridionally incompressible surfaces
that remains open is: {\it is it possible for a $(1,1)$-knot to have two
nonisotopic closed meridionally incompressible surfaces ? }  It would be natural to
conjecture a negative answer for this question. 
Such a conjecture has been verified for satellite $(1,1)$-knots in $S^3$
\cite{Sa}. In paticular this implies that if a $(1,1)$-knot in $S^3$ contains 
a closed meridionally incompressible surface of genus $\geq 2$ in its 
complement, then it is a hyperbolic knot.

We remark that there are
$(1,1)$-knots which contain closed incompressible, but meridionally compressible
surfaces. If a knot contains such a surface $S$, then by doing meridional
compressions one gets an incompressible meridional surface $\tilde S$, i.e., a
surface which intersects the knot in meridians. If enough meridional compressions
are done then the surface $\tilde S$ will be meridionally incompressible, 
that is, if
$D$ is a disk with $D\cap \tilde S =\partial D$, and which intersects $k$ in one
point, then $\partial D$ is parallel in $\tilde S$ to a boundary component of
$\tilde S$. It can be proved that for given integers $g\geq 1, \ n\geq 1$  
($g\geq 0, \ n\geq 1$), there exist
$(1,1)$-knots in $S^3$ ($S^1\times S^2$ or $L(p,q)$) which contain a meridionally
incompressible surface of genus
$g$ with $2n$ punctures \cite{ER}. A similar result for tunnel number one knots in
$S^3$ was proved previously in \cite{E3}, but those knots do not seem to be
$(1,1)$-knots. This contrasts with the main result of \cite{GR}, which shows that 
if a
tunnel number one knot contains an incompressible  meridional surface of genus
$0$ then $M$ has a $S^1\times S^2$ or $L(p,q)$ summand.

Given a tunnel number one knot $k$, it is a difficult problem to determine if $k$
is a $(1,1)$-knot. A method for determining this for knots in $S^3$ is done in \cite{MSY},
by using some computations coming from quantum invariants, and some explicit examples
are given. Another approach is given in \cite{MR}, where the existence of such knots is
proved, but no explicit example is given. We can give more examples of such
phenomena. Theorem 3 implies that if a
$(1,1)$-knot $k$ contains a closed, meridionally
incompressible surface, then the surface bounds a handlebody in $M$ in which $k$
lies. But in general, the knots constructed in Section 6 of \cite{E2} do not have 
this property.
That construction was done for knots in $S^3$, but it is not difficult to see that
it works with the obvious modifications for knots in $S^1\times S^2$ or in a lens
space $L(p,q)$. In particular we have, 

\begin{corollary} Let $k$ be a knot as in Theorem 2. Let $K$ be an
iterate of $k$ as in Theorem 6.4 of \cite{E2}. Then $K$ is a tunnel number one knot
which is not a $(1,1)$-knot.
\end{corollary}

This shows that there exists tunnel number one knots in $S^1\times S^2$ or in a lens
space $L(p,q)$ which are not $(1,1)$-knots, and explicit examples are easy to
construct. Those seem to be the first examples of that phenomenon.

\section {Proof of Theorem 3 }

Consider $T\times I \subset M$, where $T$ is a Heegaard torus for $M$, and let $T_0$,
$T_1$, $R_0$, $R_1$ be as in Section 2. Let $k$ be a $(1,1)$-knot, and assume
that $k$ lies in $T\times I$, such that $k\cap (T\times \{ 0\}) = k_0$ is an arc,
$k\cap (T\times \{ 1\}) = k_1$ is an arc, and $k\cap (T\times (0,1))$ consists of two
straight arcs. 

Suppose there is a closed meridionally incompressible
surface $S$ in $M-k$. Assume that $S$ intersects transversely $T_0$ and
$T_1$. Let $S_0=S\cap R_0$, $S_1=S\cap R_1$, and $\tilde S=S\cap (T\times I)$.

Let $\pi:T\times I \rightarrow I$ be the height function, where we choose $0$ to be 
the highest point, and $1$ the lowest. We may assume that the height function on 
$\tilde S$ is a Morse function. So there is a finite set of different points
$X=\{x_1,x_2,\dots,x_m\}$ in $I$, so that
$\tilde S$ is tangent to $T\times \{ x_i \}$ at exactly one point, and this
singularity can be a local maximum, a local minimum, or a simple saddle. For any
$y\notin X$,
$T\times \{ y \}$ intersects $\tilde S$ transversely, so for any such $y$, 
$\tilde S \cap (T\times \{ y \})$ consists of a finite collection of simple closed
curves called level curves, and at a saddle point $x$, either one level curve of
$\tilde S$ splits into two level curves, or two level curves are fused into one
curve.

Define the complexity of $S$ by the pair 
$c(S)=(\vert S_0 \vert + \vert S_1 \vert + \vert \tilde S \vert, \vert X \vert)$ 
(where $\vert Y \vert$ denotes the number of points if $Y$ is a finite set, and the
number of connected components if it is a surface; give to such pairs the
lexicographical order).
Assume that $S$ has been isotoped so that $c(S)$ is minimal.

\begin{claim} The surfaces $S_0$, $S_1$ and $\tilde S$ are
incompressible in $R_0$, $R_1$, and $(T\times I) - k$ respectively.
\end{claim}

\begin{proof}  Suppose one of the surfaces is compressible, say $\tilde S$, and let
$D$ be a compression disk, which has to be disjoint from $k$. Then $\partial D$ is
essential in $\tilde S$ but inessential in $S$. By cutting $S$ along $D$ we get a
surface $S'$ and a sphere $E$. Note that $S$ and $S'$ are isotopic in $M-k$ (this
could fail if $M=S^1 \times S^2$, but if it fails then there would be an
incompressible sphere in $M-k$, and then $k$ would be a trivial knot and $S$
a sphere). For $S'$ we can similarly define the surfaces $S_0'$, $S_1'$ and 
$\tilde S'$. Note that $\vert S_i \vert = \vert S_i' \vert + \vert E \cap R_i
\vert$, $i=1,\ 0$, so either $\vert S_0' \vert < \vert S_0 \vert$ or 
$\vert S_1' \vert < \vert S_1 \vert$, for $E$ intersects at least one of $R_0$,
$R_1$. Also $\vert \tilde S' \vert \leq \vert \tilde S \vert$, so
$c(S') < c(S)$, but this contradicts the choice of $S$.
\end{proof}

This implies that $S_0$ is a collection of trivial disks, meridian
disks and essential annuli in $R_0$. If a component of $S_0$ is a trivial disk $E$,
then $\partial E$ bounds a disk on $T_0$ which contains $k_0$, for otherwise
$\vert S_0 \vert$ could be reduced. If a component of $S_0$ is an essential annulus
$A$, then $A$ is parallel to an annulus $A^\prime \subset T_0$, and
$A^\prime$ must contain $k_0$, for otherwise $\vert S_0 \vert$ could be reduced.
This also implies that the slope of $\partial A$ cannot consist of one longitude
and several meridians of the solid torus $R_0$, for in this case
$A$ would also be parallel to $T_0-A^\prime$, and then $\vert S_0\vert$
could be reduced. It also implies that $S_0$ cannot contain both
essential annuli and meridian disks. A similar thing can be said for $S_1$.

\begin{claim} $\tilde S$ does not have any local maximum or minimum.
\end{claim}

\begin{proof} Suppose $\tilde S$ does have a maximum. It will be shown that 
$\tilde S$ has a component which is either a disk, an annulus, a once punctured
annulus, or a once punctured torus which is parallel to a subsurface in $T_1$. 
Choose the maximum at lowest level, say at level
$x_i$, so that there are no other maxima between $x_i$ and 1. So 
$\tilde S\cap (T\times \{ x_i \})$ consists of a point and a collection of simple 
closed curves. Just below the level $x_i$, the surface $\tilde S$ intersects
the level tori in simple closed curves, so below the maximum, 
a disk $E_1\subset \tilde S$ is
being formed; if $x_i$ is the last singular point then the disk $E_1$ will be
parallel to a disk on $T_1$. 

Look at the next singular level
$x_{i+1}$, $x_i<x_{i+1}< 1$. If this is a local minimum, or a saddle whose level curves
are  disjoint from the boundary curve of $E_1$, then we can interchange the
singularities.  If it is a saddle point involving the curve $\partial E_1$ and
another curve, then both singular points cancel each other. If it is a saddle formed
by an autointersection of the curve $\partial E_1$, then the disk $E_1$ below the
maximum transforms into an annulus $E_2$, which must be parallel to an annulus in
some $T\times \{y\}$, for otherwise the annulus $E_2$ would be compressible in
$T\times [0,y]$, and then either $\tilde S$ would be compressible or could be
isotoped to reduce the number of singular points. 
For any $y$ just below $ x_{i+1}$, $\partial E_2$ consists of two parallel curves
$c_1$ and $c_2$. If $c_1$ is a trivial curve on
$T\times \{ y\}$ which bounds a disk disjoint from $k$, then by cutting $\tilde S$
with such disk, we get a surface $S'$ isotopic to $S$, but with $c(S')< c(S)$, for
$S'$ has less singular points than $S$. If $c_1$ bounds a disk which intersects 
$k$ once, then $S$ would be meridionally compressible. So there remain two
possibilities, either $c_1$ and $c_2$ are essential curves on $T\times \{ y\}$ or
$c_1$ bounds a disk which intersects $k$ twice.

Look at the next singular level $x_{i+2}$. If it is a local minimum or a saddle whose
level curves are disjoint from the annulus $E_2$, then again we can interchange the
singularities (note that these curves cannot be inside the solid torus bounded by
$E_2$ and an annulus in  $T\times \{ x_{i+2} \}$, for there would be a
maximum in there). If it is a saddle joining the annulus $E_2$ with another curve,
then by pushing $E_2$ down the singular point $x_{i+2}$ is eliminated. So this
singular point has to be of the annulus
with itself, and then a new surface $E_3$ is formed. The surface $E_3$ has to be
a once punctured torus or a once punctured annulus, which is parallel to some 
$T\times\{y^\prime\}$, for otherwise $\tilde S$ would be compressible or could be
isotoped to reduce the number of singular points. If there is one more singular point
of $E_3$ with itself, then the surface
$E_4$ which will be formed would be compressible, for a component of its
boundary would bound a disk in some $T\times \{y''\}$ disjoint from $k$.  We
conclude that either
$\tilde S$ is compressible, the number of singular points is not minimal, or a
component of $\tilde S$, say
$S^\prime$, it is parallel to a surface in $T_1$. In the latter case the
surface $S^\prime$ is disjoint from $k_1$, and then it can be pushed to $R_1$,
reducing $c(S)$.
\end{proof}

At a nonsingular level $y$, $\tilde S\cap (T\times \{ y \})$ consists of simple closed
curves. If such a curve $\gamma$ is trivial in $T\times \{ y \}$, then it bounds a
disk in that torus.  If such a disk is disjoint from $k$, then by
the incompressibility of $\tilde S$, $\gamma$ bounds a disk in $\tilde S$, which give
rise to a single maximum or minimum, contradicting Claim 2. Such a disk cannot
intersect $k$ once, for
$S$ is meridionally incompressible, so $\gamma$ bounds a disk which intersects $k$
twice. 

\bigskip
\bigskip
\bigskip

\begin{claim} Only the following types of saddle points are possible.
\begin{enumerate}
\item A saddle changing a trivial simple closed curve into two essential simple
closed curves.
\item A saddle changing two parallel essential curves into a trivial curve.
\end{enumerate}
\end{claim}

\begin{proof} At a saddle, either one level curve of $\tilde S$ splits
into two level curves, or two level curves are joined into one level curve.
If a level curve is trivial in the corresponding level torus and at a saddle
the curve joins with itself, then the result must be two essential simple
closed curves, for if the curves obtained are trivial, then $\tilde S$ would be
compressible or meridionally compressible. If in a level there are nontrivial 
curves of intersection,
then there is an even number of them, for $S$ is separating. So if a curve is 
nontrivial and at the saddle joins with itself, then the result is a curve with the 
same slope as the
original and a trivial curve, for the saddle must join points on the same side of
the curve. If two trivial level curves are joined into one, then
because the curves must be concentric, the surface $S$ would be compressible. 
So only the following types of saddle points are possible:

\begin{enumerate}
\item A saddle changing a trivial simple closed curve into two essential simple
closed curves.
\item A saddle changing two parallel essential curves into a trivial curve.
\item A saddle changing an essential curve $\gamma$ into a curve with the same slope
as $\gamma$, and a trivial curve.
\item A saddle changing an essential curve $\gamma$ and a trivial curve into an
essential curve with the same slope as $\gamma$.
\end{enumerate}

We want to show that only saddles of types 1 and 2 are possible. Look at two
consecutive saddle points. We note that in many cases $\tilde S$ can be isotoped
so that two consecutive singular
points can be put in the same level, and in that case a compression disk for  
$\tilde S$ can be found. Namely,

\begin{enumerate}
\item A type 3 is followed by a type 1, see Figure 3(a).
\item A type 3 followed by a type 2, see Figure 3(b).
\item A type 1 is followed by a type 4, see Figure 3(c).
\item A type 2 followed by a type 4, see Figure 3(d).
\item A type 3 followed by a type 4, see Figure 3(e) and 3(f).
\end{enumerate}

\bigskip

\centerline{\includegraphics{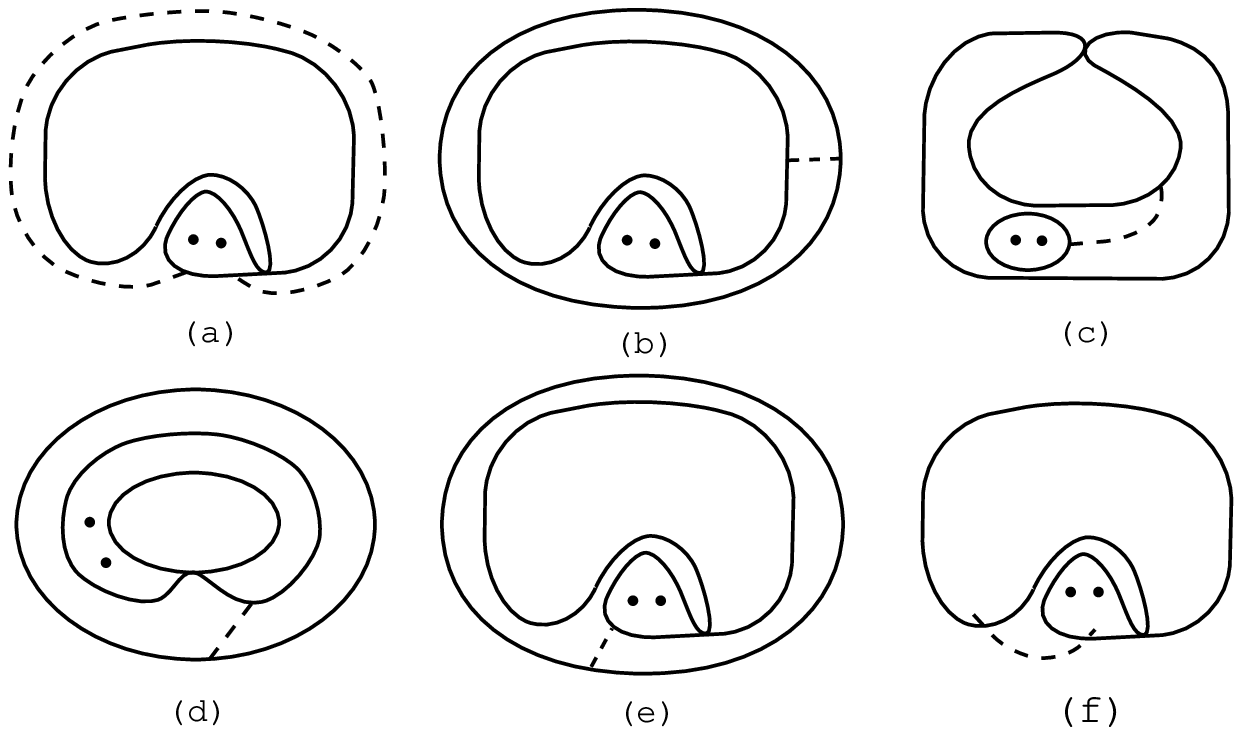}}

\centerline{\sr Figure 3}

\bigskip

The dots in Figure 3 indicate the points of intersection between $k$ and some 
$T\times \{ y \}$, and the dotted line indicates the path followed by the second
saddle point. Note that in all the other cases it may not be possible to put the
two saddles points in the same level. 

Suppose there are singularities of type 3, and take the lowest one, say at level
$x_m$. If there is a singularity at level $x_{m+1}$, then it is of type 1, 2, or 4,
which contradicts the previous observation. So $x_m$ is in fact the last
singularity.  Then $S_1$ consists of at least a disk, and at least one annulus or a
meridian disk. Isotope $S$ so that the saddle at $x_m$ is pushed into $R_1$, this
can be done without moving $k_1$. If
$S_1$ contained a meridian disk and a trivial disk, which joined in $x_m$, then now
there is just one meridian disk, reducing $c(S)$ by 2. If $S_1$ contained an annulus
and a disk, which joined in the singularity, then now there is just one annulus,
reducing $c(S)$ by 2.

If there is a saddle of type 4, then take the highest one. If it is not a level
$x_1$, then it is preceded by a singularity of type 1, 2, or 3, which is
a contradiction. So it is a level $x_1$. Push the singularity to
$R_0$. Doing an argument as in the previous case we get a contradiction.
\end{proof}

\begin{claim} $S_0$ ($S_1$) does not contain any meridian disk
of $R_0$ ($R_1$).
\end{claim}

\begin{proof} Suppose $S_0$ contains meridian disks, and possibly some trivial disks.
If $\tilde S$ has no singular points, then it would be a collection of annuli, and
$S$ would be a sphere; suppose then that there are some singularities on
$\tilde S$. The first saddle is of type 1 or 2, as in Claim 3. If it is of type 1,
then a trivial disk joins to itself. This singularity can be pushed to
$S_0$, so now there is an annulus in $S_0$ of meridional slope, but this shows that
$S$ is compressible. If the first singularity is of type 2, then two of the meridian
disks are joined by a saddle. Push this singularity to
$R_0$, possibly moving the arc $k_0$. Then the two disks are transformed into a
trivial disk. The arc $k_0$ is in the ball bounded by such trivial disk and a disk in
$T_0$, so $k_0$ can be rearranged to be level. This reduces $\vert S_0 \vert$ by one. 
If $S_1$ contains meridian disks a similar argument can be done.
\end{proof}

\begin{claim} $S_0$ ($S_1$) consists only of annuli.
\end{claim}

\begin{proof} Suppose $S_0$ contains a trivial disk. Look at the first singularity of 
$\tilde S$.
If it is of type 2, then there is an annulus in $S_0$, and its boundary components
are joined in the saddle. So just after the saddle there is a disk $E$ contained in a 
level torus whose boundary is a nontrivial curve on $\tilde S$; this disk intersects $k$, 
but because $S_0$ contains a trivial disk, $E$ can be pushed above this trivial disk 
to avoid intersections with $S$ or $k$, so it would be a compression disk for $S$. 
See Figure 4(a). So the first singularity has to be of
type 1, and a trivial disk touches with itself. Push the singularity to $R_0$, so
instead of the disk we have now an annulus, this leaves $\vert S_0\vert$ unchanged,
but reduces $\vert X\vert$ by one.
\end{proof}

\bigskip

\centerline{\includegraphics{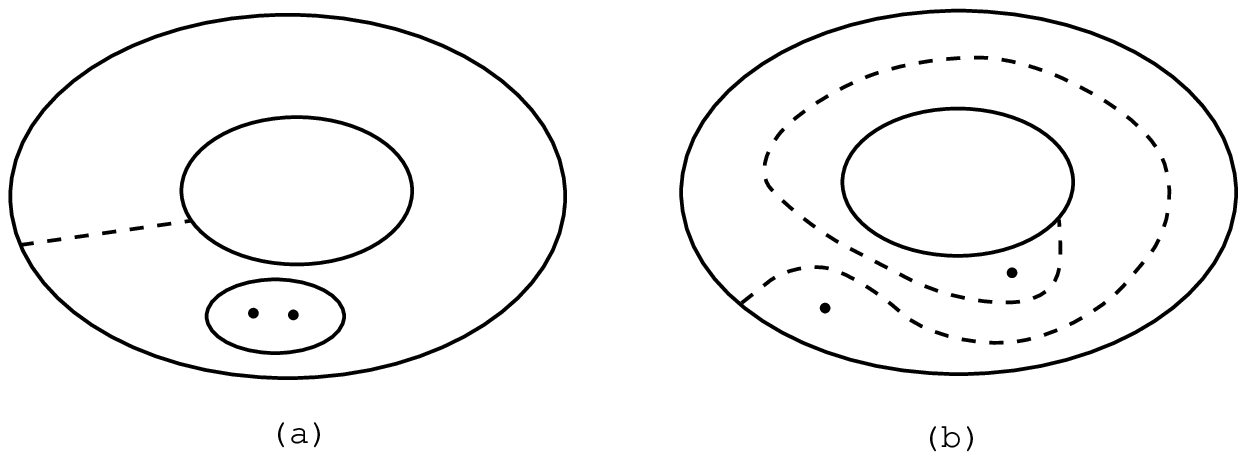}}

\centerline{\sr Figure 4}

\bigskip

Note that if the first singularity is of type 2, and $S_0$ does not contain any
disk, then the path followed by the singularity may be quite complicated, as shown
in Figure 4(b), and then a disk as the one found in Figure 4(a) is not a compression
disk for $S$, for it intersects $k$.

\begin{claim} $S_0$ ($S_1$) is just one annulus.
\end{claim}
\begin{proof} Suppose $S_0$ consists of $p$ annuli, $p>1$. These annuli are nested, 
with the innermost one, say $S_0'$,  cobounding with an annulus in $T_0$ a solid torus 
with interior disjoint from $S$ and containing an arc of $k$. As $S_0$ consists only of 
annuli, the first saddle is of type 2, and it has to join the two boundary curves of the 
innermost annuli $S_0'$. The next singularity must be also of type 2,
joining now the curves of the second innermost annuli, for if the next singularity is of 
type 1, then it is not difficult to see that there is a compression disk for $\tilde S$.
So the first $p$ singularities are of type 2.
At a level just after the first $p$ singularities, the intersection of $\tilde S$ with that 
level torus consists of $p$ nested trivial curves, the innermost one bounding a disk which 
intersects $k$ in two points. The next singularity must be of type 1,
and joins the outermost of such curves with itself. The next $p-1$ singularities are also 
of type 1, for otherwise there would be a compression disk for $\tilde S$.
At a level just after these singularities, the intersection of $S$ with that level 
consists of $2p$ essential curves. Two of these curves bound an annulus in the level 
torus which intersects $k$ in two points; note that these curves are in the same 
component of $\tilde S$ as the boundary curves of $S_0'$.

This shows that if $S_0$ consists of $p$ annuli,
then the first $p$ singularities are of type 2, which are followed by $p$
singularities of type 1, and then are followed by $p$ singularities of type 2, and so on. 
This implies in particular that $S_0 \cup \tilde S$ consists of $p$ 
connected components. Note that the component which contains the annulus $S_0'$ remains 
always being innermost, that is, at each level torus, its curves of intersection with that 
torus bound an annulus or a disk which intersects $k$ in two points.
Then in $T_1$, this component bound an annulus which contains an arc of $k$.
This implies that in $R_1$ there is an annulus component of $S_1$, say $S_1'$, having 
these curves as its boundary components, that is, $S_1'$ is an innermost annulus in 
$R_1$, so $S_0'\cup \tilde S \cup S_1'$ is a connected component of $S$,
which contradicts the fact that $S$ were connected. So we conclude that
$S_0$ is just one annulus, similarly $S_1$. 
\end{proof}

\begin{claim} If $\vert X \vert =0$, then $S=\partial N(\Gamma)$, for a
toroidal graph $\Gamma$ of type 1, and $k$ is level with respect to $\Gamma$.
\end{claim}

\begin{proof} As $S_0$ and $S_1$ are annuli, and $\tilde S$ has no singular points,
it follows that $S$ is a torus. Then $S_0$ and $S_1$ determine the same slope in
$T_0$ and $T_1$ respectively, so $S=\partial N(\gamma_1)$, where $\gamma_1$ is a
curve in $T_0$ lying in the annulus parallel to $S_0$. It follows from the
definitions that $k$ is level with respect to $\gamma_1$.
\end{proof}

In what follows suppose that $\vert X \vert \geq 1$.
The proof of Claim 6 implies the following:

\begin{claim} $\tilde S$ has an even number of saddle points. The sequence of
saddle points can be written as
$\{y_1,z_1,y_2,z_2,\dots, y_n,z_n\}$, where $y_i$ is a saddle of
type 2, and $z_i$ is a saddle of type 1. So in a level $y_i < v < z_i$,
$\tilde S\cap T\times \{ v \}$ consists of a trivial curve on $T\times \{ v \}$,
and in levels $0\leq u < y_1$, $z_i < u < y_{i+1}$ and $z_n < u \leq 1$,
$\tilde S\cap T\times \{ u \}$ consists of two essential curves on
$T\times \{ u \}$, say of slope $r_i$.
\end{claim}

\begin{claim} $S=N(\Gamma)$ for a certain toroidal graph $\Gamma$,
and $k$ can be isotoped to be level with respect to $\Gamma$.
\end{claim}

\begin{proof} Choose nonsingular levels $u_1,u_2,\dots,u_{n+1}$, so that 
$u_1=0$, $z_i < u_{i+1} < y_{i+1}$, and $u_{n+1}=1$. So 
$\tilde S\cap (T\times \{ u_i\})$ consists of two simple closed essential
curves of slope $r_i$. These curves bound an annulus in
$T\times \{ u_i \}$ which intersects $k$ in two points. Let $\gamma_i$ be an
essential curve in the interior of this annulus. 

Choose nonsingular levels $v_1,\dots,v_n$, so that $y_i < v_i < z_i$. 
So $\tilde S\cap (T\times \{ v_i \})$ consists of a simple closed curve which 
bounds a disk $D_i$ in $T\times \{ v_i \}$ which intersects $k$ in two points. 

Let $\alpha_i$ be an straight arc in $[u_i,u_{i+1}]$,
joining $\gamma_i$ and $\gamma_{i+1}$, and which is disjoint from $\tilde S$.
Let $\Gamma$ be the 1-complex consisting of the union of all the
curves $\gamma_i$ and the arcs $\alpha_j$. So $\Gamma$ is a trivalent
graph embedded in $M$. 
Do a level preserving ambient isotopy in $T\times I$,
keeping $\Gamma$ fixed and disjoint form $\tilde S$, so that at a level $v_i$, the
disk $D_i$ is isotoped to a tiny neighborhood of the point 
$D_i\cap \alpha_i$, and at a level $u_i$, the two curves in
$\tilde S\cap (T\times \{ u_i \})$ bound an annulus which is a tiny neighborhood of
the curve $\gamma_i$. 
Between the levels $v_i$ and $v_{i+1}$, the surface $\tilde S$ and the disks $D_i$
and $D_{i+1}$ bound a solid torus $P_i$, which is foliated by disks, annuli and 
two singular disks; it is foliated by disks between the levels $v_i$ and $z_i$,
and between $y_{i+1}$ and $v_{i+1}$, it is foliated by annuli between the levels 
$z_i$ and $y_{i+1}$, and there are singular disks at the levels $z_i$ and $y_{i+1}$. At the 
top of $P_i$ is the disk $D_i$, at the bottom the disk $D_{i+1}$, and at a middle
level the annulus $\eta(\gamma_i)$.
It follows that $S$ is in fact the boundary of a regular
neighborhood of $\Gamma$. The slopes of the curves $\gamma_i,\ \gamma_{i+1}$ must
satisfy $\Delta(r_i,r_{i+1})\geq 2$, for otherwise the surface $S$ would be
compressible. Then $\Gamma$ is a toroidal graph of type $n+1$.

When doing the isotopy the knot $k$ is also modified, but because the isotopy is
level preserving, the arcs $k\cap (T\times I)$ are still straight arcs.
Note that $k\cap P_i$ consists of two straight arcs with
endpoints in $D_i$ and $D_{i+1}$, so each of these arcs intersects each surface in
the foliation of $P_i$ exactly once. So $k\cap P_i$ is like a 2-braid in $P_i$.
$P_i\cap \alpha_{i-1}$ is a straight arc with one endpoint in $D_i$ and the other
in $\gamma_i$, and $P_i\cap \alpha_i$ is a straight arc with one endpoint in
$D_{i+1}$ and the other in $\gamma_i$. Now $k\cap P_i$ can be isotoped, keeping its
endpoints in $D_i$ and $D_{i+1}$, so that $k\cap N(\alpha_{i-1})$ and 
$k\cap N(\alpha_i)$
both consist of two straight arcs, and $k\cap \eta(\gamma_i)$ consists of two arcs.

Above the level $v_1$, the annulus $S_0$, the surface $\tilde S$ and the disk
$D_1$ bound a solid torus $P_0$. $k\cap P_0$ consists of two straight arcs going
from $D_1$ to $\eta(\gamma_0)$, and an arc lying in $\eta (\gamma_0)$.
$P_0 \cap \alpha_1$ consist of a straight arc going from $D_1$ to $\gamma_0$.
Now it is not difficult to see that $k\cap P_0$ can be isotoped so that it
consists of two straight arcs lying in $N(\alpha_1)$, and an arc in $\eta (\gamma_0)$.
A similar argument can be done for $P_n$, the solid torus lying below the disk
$D_n$. This shows that $k$ is level with respect to $\Gamma$.
\end{proof}

The proof of Theorem 3 is now complete.

\bigskip

{\br{Acknowledgement}.}
Part of this work was done while visiting the Department of Mathematics of The
University of California at Santa Barbara in 1999. I am grateful for their
hospitality. This research was partially supported by a NSF grant, and by a
sabbatical fellowship from CONACYT. I am also grateful to the anonymous 
referee for
doing valuable suggestions.

\end{document}